\pdfoutput=1

\documentclass[11pt]{article}
\usepackage{a4}
\usepackage{amsmath,amsfonts,amssymb,amscd,amsthm}
\usepackage[pdftex]{graphicx}
\usepackage{subfigure}
\usepackage{pict2e}
\usepackage{color}
\usepackage{hyperref}

\newtheorem{theorem}{Theorem}

\newtheorem{example}{Example}

\newtheorem{alg}{Algorithm}

\newcommand{\Keywords}[1]{\par\noindent
{\small{ \textbf{Key words}\/}: #1}}

\begin{document}
\title{Numerical Simulations for the Non-linear Molodensky Problem}
\author{L.~Banz \thanks{Institut f\"ur Angewandte Mathematik,
           Leibniz Universit\"at Hannover,
           Welfengarten 1, 30167 Hannover, Germany.
           Email: {\tt costea@ifam.uni-hannover.de},
									{\tt banz@ifam.uni-hannover.de},
                  {\tt stephan@ifam.uni-hannover.de}}
\and A.~Costea$^*$
\and H.~Gimperlein
   \thanks{Department of Mathematical Sciences, University of Copenhagen, Universitetsparken 5, 2100 Copenhagen \O, Denmark.
    Email:
    {\tt gimperlein@math.ku.dk}}
\and E.P.~Stephan$^*$}

\maketitle
\begin{abstract}
We present a boundary element method to compute numerical approximations to the non-linear Molodensky problem, which reconstructs the surface of the earth from the gravitational potential and the gravity vector. Our solution procedure solves a sequence of exterior oblique Robin problems and is based on a Nash-H{\"o}rmander iteration. We apply smoothing with the heat equation to overcome a loss of derivatives in the surface update. Numerical results compare the error between the approximation and the exact solution in a model problem.\\
\end{abstract}

\section{Introduction}
The determination of the shape of the earth and its exterior gravitational field from terrestrial measurements is a basic problem in physical geodesy \cite{c3,Heiskanen1967}. Molodensky \cite{Molodensky1,Molodensky2} formulated it as an exterior free boundary problem for the Laplace equation in $\mathbb{R}^3$ with the gravitational potential $W$ and field $G$ prescribed on an unknown boundary diffeomorphic to the two--dimensional sphere by a map $\varphi: S^2 \rightarrow \mathbb{R}^3$. With the advent of satellite technologies to determine the surface of the earth high-precision studies combine satellite data with local gravity measurements.\\
H{\"o}rmander \cite{Hoermander} proved local existence and uniqueness of the solution of Molodensky's problem. Based on ideas of Nash and Moser, his constructive proof overcomes the loss of regularity of subsequent iterates in standard fixed-point methods for this problem by introducing an additional smoothing operator in each step. In \cite{Costea2013} we have shown that smoothing by a higher-order heat equation can be used and is numerical feasible. We have obtained error estimates for the iterates showing the dependence of the rate of convergence of the algorithm on certain parameters.\\
Here we present some computational aspects of our approach to the non-linear Molodensky problem. To solve the free boundary problem we iteratively construct a sequence $(\varphi_m)_{m\in \mathbb{N}_0}$ of approximations to the boundary $\varphi$, where $\varphi_m$ is obtained from the boundary element solution of the problem linearized around $\varphi_{m-1}$. The numerical solution of the linearized Molodensky problem using the boundary element method with piecewise linear ansatz functions was first analyzed by Klees et al.~\cite{Schwabfast}, see also \cite{Holota,FreedenMayer06}. Smoothed iterative solvers involving the heat equation were investigated in \cite{Jerome,Jeromeheatsm} for ordinary differential equations.\\
In Section~\ref{sec:2} we present the Nash-H\"ormander algorithm and in Section~\ref{sec:3} the boundary element procedure. Section~\ref{numex} gives the detailed description for a model problem in which the sphere of radius 1.1 is recovered by our method starting from the unit sphere. \\
However, our solution procedure can be directly applied to more complicated geometries like spheroids, telluroids, etc. In particular, to be relevant for the geodetic community, one might model the exact surface by the ETOPO1 model of the earth and compute the gravity vector from the EGM2008 model. A more realistic model problem could then try to recover $\varphi$ starting from the GRS80 ellipsoid as surface $\varphi_0$ and the corresponding Somigliani-Pizetti field as $(W_0,G_0)$ \cite{Ardalan2001}.

\section {The Nash-H\"ormander algorithm}
\label{sec:2}
\hspace*{5mm} A classical problem in geodesy is to find an embedding $\varphi: \mathbb S^2 \rightarrow  \mathbb R^3$ such that $G=\Gamma(W,\varphi)$ where the potential $W$ and the gravity vector $G$ are given on $ \mathbb S^2$. The nonlinear map $\Gamma$ is implicitly described by the nonlinear Molodensky problem \eqref{eq:NLmol}:
 Find $\varphi:\mathbb S^2 \rightarrow \mathbb R^3$ s.t.
\begin{align}
\Delta w&=0 \quad \mbox {outside} \quad \varphi(\mathbb S^2) \nonumber \\
\quad w &=W \circ \varphi^{-1},\ g:=\nabla w = G \circ \varphi^{-1}\quad \mbox {on} \quad  \varphi(\mathbb S^2) \label{eq:NLmol}\\
 w(x) &= \frac{M}{|x|} +O(|x|^{-3}) \quad \mbox{when} \quad |x| \rightarrow \infty,\quad M \in \mathbb R. \nonumber
\end{align}
In this model the rotation of the earth is neglected, and the decay condition fixes the center of mass
of the earth $\varphi(\mathbb S^2)$ at the origin.\\
\hspace*{5mm} In \cite{Hoermander} the nonlinear Molodensky problem \eqref{eq:NLmol} is solved by a sequence of linearized problems together with a specific smoothing of a Newton--like iteration. H{\"o}rmander shows convergence of his sequence of approximate solutions. A rate of convergence of the iteration was derived in \cite{Costea2013}, and a boundary element procedure was suggested to solve \eqref{eq:NLmol}. Here we present numerical experiments with the solution procedure from \cite{Costea2013} and describe how to apply the heat equation as a smoother within this context. In each step of the iteration one solves the linearized Molodensky problem:
Find $u$ and $a_{j} \in \mathbb R$ such that
\begin{align}
&\Delta u=0 \quad \mbox{outside} \quad \varphi(S^2), \nonumber \\
& u + \nabla u \cdot h = F(G,W) -\sum_{j=1 }^3 a_{j} A_j  \quad \mbox{on} \quad \varphi(S^2)   \label{eq:linNLmol}\\
& u(x)=\frac{c}{|x|} + O(|x|^{-3}) \quad \mbox{when} \quad |x| \rightarrow \infty,\quad c \in \mathbb R , \nonumber
\end{align}
where the Marussi condition $\det(\nabla g) \neq 0$ on $\varphi(S^2)$ for  $g=\nabla U$ with a fixed  $U$ is assumed to hold. The vector $h$ in \eqref{eq:linNLmol} satisfies
$h = -(\nabla g)^{-1} g$. The functions $A_j = \frac{x_j}{|x|^3}$ guarantee that \eqref{eq:linNLmol} is well-posed under certain conditions, see \cite{Hoermander,Costea2013}.\\
Then one updates $\varphi$ by the increment $\Delta \dot{\varphi}=\Delta (\nabla g \circ \varphi)^{-1}(\dot{G} - \nabla u \circ \varphi)$, where $\dot{G}$ is the increment of $G$. One observes as main difficulty in constructing $\varphi$ that $\dot{\varphi}$ is less regular than  $\varphi$. We overcome this dilemma by applying the heat-equation as described below obtaining the smoothed quantities $\dot{\widetilde W}_m$ and $\dot{\widetilde G}_m$. The Nash-H\"ormander iteration for \eqref{eq:NLmol} reads as follows:
At each iteration step, for $F_m = \dot {\widetilde W}_m\circ \varphi_m^{-1}+(\dot {\widetilde G}_m \circ \varphi_m^{-1}) \cdot h_m$ (as specified below) find $u_{m}$ and $a_{j,m} \in \mathbb R$ such that
\begin{align}
&\Delta u_{m}=0 \quad \mbox{outside} \quad \varphi_{m}(S^2),  \nonumber \\
& u_{m} + \nabla u_{m} \cdot h_{m} = F_{m}(G,W) -\sum_{j=1 }^3 a_{j,m} A_j  \quad \mbox{on} \quad \varphi_{m}(S^2)  \label{eq:h_equation} \\
& u_{m}(x)=\frac{c}{|x|} + O(|x|^{-3}) \quad \mbox{when} \quad |x| \rightarrow \infty,\quad c \in \mathbb R. \nonumber
\end{align}
with
 $\det(\nabla g_{m-1}) \neq 0$ on $\varphi_{m}(S^2)$ and $g_{m-1}=\nabla w_{m-1} = \nabla \sum_{j=0}^{m-1}u_{{j}}$
and
\begin{align*}
h_{m} = (-(\nabla g_{m-1})^{-1} g_{m-1} \circ \varphi_{m-1})\circ \varphi_{m}^{-1},\quad  A_j = \frac{x_j}{|x|^3}.
\end{align*}
The solution $u_m$ determines the nonlinear correction $\dot{\varphi}_m$, i.e.~$\varphi_{m+1}=\varphi_m +\delta_m \dot{\varphi}_m$, namely
\begin{align}
\dot{\varphi}_{m}=(\nabla g_{m} \circ \varphi_{m})^{-1}(\dot{\widetilde G}_{m} - \nabla u_{m} \circ \varphi_{m}). \label{corrphi}
\end{align}
The gravitational vector $g_m$ is determined from the approximation to the potential as computed in the first $m$ steps,
 \begin{equation}
 w_m=W_{m-1} \circ \varphi_m^{-1} +\triangle_m u_m \quad \mbox{on} \quad \varphi_m(S^2)\ , \label{wmtotal}
\end{equation}
$$W_{m-1}=\begin{cases}v_0 \circ \varphi_0 + \triangle_0 u_0 \circ \varphi_0 + \triangle_1 u_1 \circ \varphi_1 + \dots + \triangle_{m-1} u_{m-1} \circ \varphi_{m-1},\! \quad\mbox {for} \quad\!m\geq 1\\
v_0 \circ \varphi_0  \quad\mbox {for} \quad m=0 \end{cases}  \ ,$$
for suitable stepsizes $\Delta_j$ and initial approximation $v_0$ by solving an exterior Dirichlet problem:\\
For given $w_m$ on $\varphi_m(S^2)$, find $ \overline {v}_m : \mathbb R^3\backslash \bar \Omega_m \rightarrow \mathbb R$ and constants $a_{j,m} \in \mathbb R$ such that
\begin{align}
&\Delta \overline {v}_m =0 \quad \mbox{in} \quad \mathbb R^3\backslash \bar \Omega_m, \nonumber\\
& {\overline{v}_m}_{|_{\partial {\Omega_m}}} = w_m -\sum_{j=1}^3 a_{j,m}  A_j (x)\big|_{x \in\:\varphi_m(S^2)}  \quad \mbox {on} \quad \varphi_m(S^2) \label{Dirichlet} \\
& \overline{v}_m(x)= \frac {c}{|x|} + O(|x|^{-3}) \quad \mbox {when} \quad  |x| \rightarrow \infty \ .\nonumber
\end{align}
Now with (\ref{corrphi}) the surface update $\dot \varphi_m$ is computed using $g_m=\nabla \overline{v}_m$ and $\nabla g_m=\nabla^2 \overline{v}_m$.\\
The full iterative method involves smoothing in each step based on the solution operator $S_{\theta}$ to a higher-order heat equation (discussed below). It reads as follows:
\begin{alg} \label{alg:nashhormsm}
\vspace*{10mm} \hrule
\vspace*{0.5mm} (Nash-H\"ormander algorithm) \vspace*{2mm}
\hrule
\begin{enumerate}
	\item For given measured data $W, G$, choose $W_0, G_0,h_0,\varphi_0, \theta_0 \gg 1, \kappa \gg 1$
  \item For $m=0,1,2,\ldots$ do
  \begin{enumerate}
        \item Compute \begin{equation}\label{thetatocompute}
                                      \theta_m=(\theta_0^{\kappa}+m)^{1/\kappa},\quad \triangle_m=\theta_{m+1}-\theta_{m}
                                    \end{equation}
         \item Compute
                 \begin{align}\label{wreccurencetosmooth}
                  \dot {\widetilde W}_0:&=S_{\theta_0}\dot W_0=S_{\theta_0}\big(\frac{W-W_0}{\triangle_0} \big) \nonumber\\
                  \dot {\widetilde W}_m:&=\frac{1}{\triangle_m}\big(S_{\theta_m}\!(W\!-W_0\! )\!-\!S_{\theta_{m-1}}\!(W \!- W_{0})\big)
                 \end{align}
        \item Compute
                 \begin{align}\label{greccurencetosmooth}
                  \dot {\widetilde G}_0:&=S_{\theta_0}\dot G_0=S_{\theta_0}\big(\frac{G-G_0}{\triangle_0} \big) \nonumber\\
                  \dot {\widetilde G}_m:&=\frac{1}{\triangle_m}\big(S_{\theta_m}\!(G\!-G_m\!+\!\! \sum_{j=0}^{m-1}\!\!\triangle_j \dot {\widetilde G}_j )\!-\!S_{\theta_{m-1}}\!(G \!- G_{m-1} +\!\! \sum_{j=0}^{m-2}\!\!\triangle_j \dot {\widetilde G}_j)\big)
                 \end{align}
        \item Find $u_m$ by solving the linearized problem \eqref{eq:h_equation} with $(\dot W_m, \dot G_m)$ replaced by $(\dot {\widetilde W}_m, \dot{\widetilde G}_m)$
        \item Find $\overline v_m$ by solving (\ref{Dirichlet}) with $w_m$ as defined in  (\ref{wmtotal})
        \item Compute $g_m=\nabla \overline {v}_m$ and $\nabla g_m=\nabla^2 \overline {v}_m$
        \item Compute the surface increment $\dot \varphi_m$ by $$\dot \varphi_m= (\nabla g_m \circ \varphi_m)^{-1} (\dot{\widetilde G}_m- \nabla u_m \circ \varphi_m)$$
                  and update surface map by $\varphi_{m+1}=\varphi_m + \triangle_m \dot \varphi_m$
        \item Update direction vector and gravity potential by
     \begin{align*}
      h_{m+1}&=((-(\nabla g_m)^{-1} g_m) \circ \varphi_m) \circ ({\varphi_{m+1}})^{-1} \\
      G_{m+1}&=g_m \circ \varphi_m
     \end{align*}
      \item Stop if $\left\| g_m\circ \varphi_m-G \right\| +\left\| \overline v_m\circ \varphi_m-W \right\|< \text{tol}$
    \end{enumerate}
\end{enumerate}
\hrule
\end{alg}
$\|\cdot\|$ might usually be chosen to be e.g.~an $\mathcal{H}^a$--norm.\\

In \cite{Costea2013} the following convergence result for the Nash-H\"ormander iteration (Alg~\ref{alg:nashhormsm}) is proved, where it is assumed that the starting values $W_0$ and $G_0$ are already in a small neighborhood (in the H{\"o}lder space $C^{\alpha+\epsilon}$) of the final values $W$ and $G$.\\
\begin{theorem}
For $\alpha>2+2\epsilon$, $0<a <\alpha$ and $\tau>0$ small, such that $a-\alpha-\tau < 0$ there exist constants $\theta_0>0$, $C_\tau>0$ s.~t.~$\varphi_m$ satisfy for all $m \geq 0$
\begin{equation*}
\|\varphi-\varphi_m\|_{C^{a+\epsilon}}\leq C_\tau \left(\|W-W_0\|_{C^{\alpha+\epsilon}}+\|G-G_0\|_{C^{\alpha+\epsilon}}\right)\theta_m^{a-\alpha+\tau}
\end{equation*}
\end{theorem}

\section {Boundary element procedure}
\label{sec:3}
Inserting a single layer potential ansatz
\begin{align}
u(x)=V \mu(x) = -\frac{1}{4 \pi}\int_{\varphi_m(\mathbb S^2)} \frac{\mu(y)}{|x-y|} ds_y \label{eq:SingleAnsatz}
\end{align}
into the linearized Molodensky problem \eqref{eq:linNLmol} translates the oblique Robin boundary condition to a second kind integral equation for the density $\mu$ on the surface $\varphi(\mathbb S^2)$, namely
\begin{align}
B\mu(x):=\frac{1}{2} \cos \beta \mu(x) +K'(\boldsymbol h)\mu(x)+V \mu(x)=f(x), \label{definitionofS}
\end{align}
where $\beta =\measuredangle(\boldsymbol n,\boldsymbol h),  K'(\boldsymbol h) = h \cdot \nabla V$ and $f=F(G,W)-\sum_{j=1}^3 a_j A_j$.\\
The following variational formulation of the above integral equation is solved by boundary elements:
Find $(\mu_{m,h},a_{j,m}^h) \in S_{h,m} \times \mathbb R^3$ s.t.
\begin{alignat}{2}
\langle B \mu_{m,h}, \psi_h \rangle_{\varphi_m^h(\mathbb S^2)}  + \sum_{j=1}^3\langle B a_{j,m}^h A_j , \psi_h \rangle_{\varphi_m^h(\mathbb S^2)} &=\langle F_m^h,\psi_h \rangle_{\varphi_m^h(\mathbb S^2)}& &\,\,\,\forall\: \psi_h \in S_{h,m}\nonumber \\
\langle \mu_{m,h},A_k \rangle_{\varphi_m^h(\mathbb S^2)} &=0&  &\,\,\,k=1,2,3 \label{eq:discreteform1}
\end{alignat}
where $S_{h,m}$ denotes the set of piecewise quadratic continuous functions on the approximate surface
$\varphi_m^h(\mathbb S^2)$ and $F_m^h$ is the projection of $F_m$ onto $\varphi_m^h(\mathbb S^2)$. This surface is obtained starting from an initial regular mesh of triangles.\\
The auxiliary exterior Dirichlet problem \eqref{Dirichlet} is again solved by the boundary element method with a single layer potential ansatz:

Find $(\widetilde \mu_{m,h},\widetilde a_{j,m}^h) \in S_{h,m} \times \mathbb R^3$ s.t.
\begin{alignat}{2}
\langle V \widetilde \mu_{m,h}, \xi_h \rangle_{\varphi_m^h(\mathbb S^2)}  + \sum_{j=1}^3
\widetilde a_{j,m}^h
 \langle V  A_j , \xi_h \rangle_{\varphi_m^h(\mathbb S^2)} &=\langle w_m^h,\xi_h \rangle_{\varphi_m^h(\mathbb S^2)}& &\,\,\,\forall\: \psi_h \in S_{h,m}\nonumber \\
\langle \widetilde \mu_{m,h},\widetilde A_k \rangle_{\varphi_m^h(\mathbb S^2)} &=0&  &\,\,\,k=1,2,3 \label{eq:discreteform2}
\end{alignat}
where $w_m^h$ is the projection of $w_m$ on $\varphi_m^h(\mathbb S^2)$. The convergence as $h \rightarrow 0$ of the boundary element approximations $(\mu_{m,h},a_{j,m}^h)$ and $(\widetilde \mu_{m,h},\widetilde a_{j,m}^h)$ is shown in \cite{Schwabfast}.\\
The crucial point for the numerical computation of the updates in the Nash-H{\"o}rmander algorithm is the Hessian $\nabla^2 \bar{v}_m$ on the surface $\varphi_m^h(\mathbb S^2)$, which must be approximated very accurately. This is difficult as the single layer potential ansatz for $\bar{v}_m$ leads to the evaluation of hypersingular integrals
\begin{equation*}
\nabla^2 \bar{v}_m=-\frac{1}{4\pi} p.f.~\int_{\varphi_m(\mathbb S^2)} \!\!\nabla_x^2 \frac{1}{|x-y|} \bar{\mu}_m (y) ds_y, \!\!\quad \! x\in \varphi_m(\mathbb S^2).
\end{equation*}
Since the gradient $g=\nabla \bar{v}=\nabla V \overline \mu_m$ can be computed analytically on plane surface pieces (see \cite{MaischakTech}), we approximate the second derivative of $\bar{v}$ by appropriate
finite differences. For the normal and tangential derivative of $g$ we take
\begin{align*}
\frac{\partial g(x)}{\partial n}&=\frac{4g(x + \delta \cdot n) - 3g(x) - g(x + 2\delta \cdot n) }{2\delta} + O(\delta^2) \\
\frac{\partial g(x)}{\partial t} &=\frac{g(x+\delta \cdot t) -g(x-\delta \cdot t) }{2\delta}  + O(\delta^2)
\end{align*}
with stepsize $\delta$. Second and higher--order derivatives have also been analyzed by Schulz, Schwab and Wendland, e.g.~in \cite{SchulzSchwabWendland,SchwabWendland}. They compute the second normal derivative from the less singular tangential derivatives using geometrically graded meshes near the singularity. For the current problem, a simpler approach gives sufficient accuracy.\\

\begin{figure}[tbp]
	\centering
	\includegraphics[keepaspectratio,width=120mm]{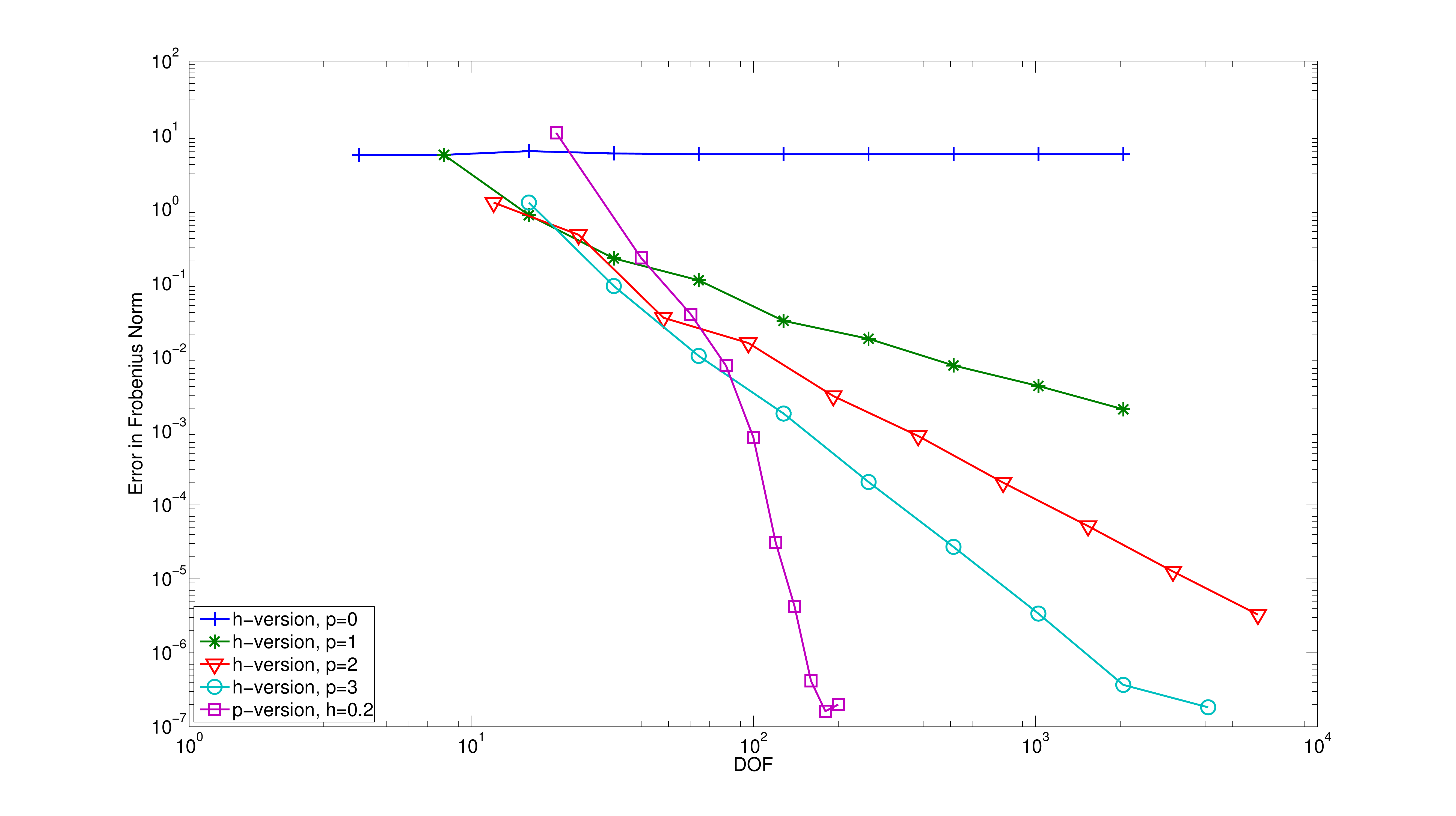}
	\caption{Error of the Hessian approximation in 2D for a point on the boundary surface}
	\label{fig:hessian2D}
\end{figure}

In the computations, the step size $\delta$ is set to $10^{-4}$ for the normal component and to $10^{-5}$ for the tangential component when approximating second derivatives.
For the presented numerical experiments the FD-approximation error is of magnitude $10^{-7}$ if no Galerkin-BEM  approximation error were to occur.
However, for very small step sizes the finite differences become numerically instable, and for the given BE-spaces the BEM-error dominates the FD-error. If $H=(H_{ij})$ denotes the exact and $H_h=(H_{h,ij})$ the approximated Hessian, we measure the error in a point $x$ as $\left(\sum_{i,j} (H_{h,ij}(x)-H_{ij}(x))^2\right)^{1/2}$.
\begin{example}
Let $\Omega=\left[-\frac{1}{2},\frac{1}{2}\right]^2$ be the domain and $u=\ln \|x\|$ the exact solution. Then the exact Hessian is $H(x)= \frac{1}{x^2} \left(
\begin{array}[pos]{c c}
	1-2x_1^2 & -2x_1x_2 \\ -2x_1x_2 & 1-2x_2^2
\end{array}\right)$.
Figure~\ref{fig:hessian2D} shows the pointwise error of the Hessian approximation in the point $x=(\frac{1}{2},\frac{1}{3})$ for $h$--versions of BEM with polynomial degree $p=0,1,2,3$, as well as for a $p$--version with $h=0.2$.
\end{example}

\begin{figure}[tbp]
	\centering
	\includegraphics[keepaspectratio,width=120mm]{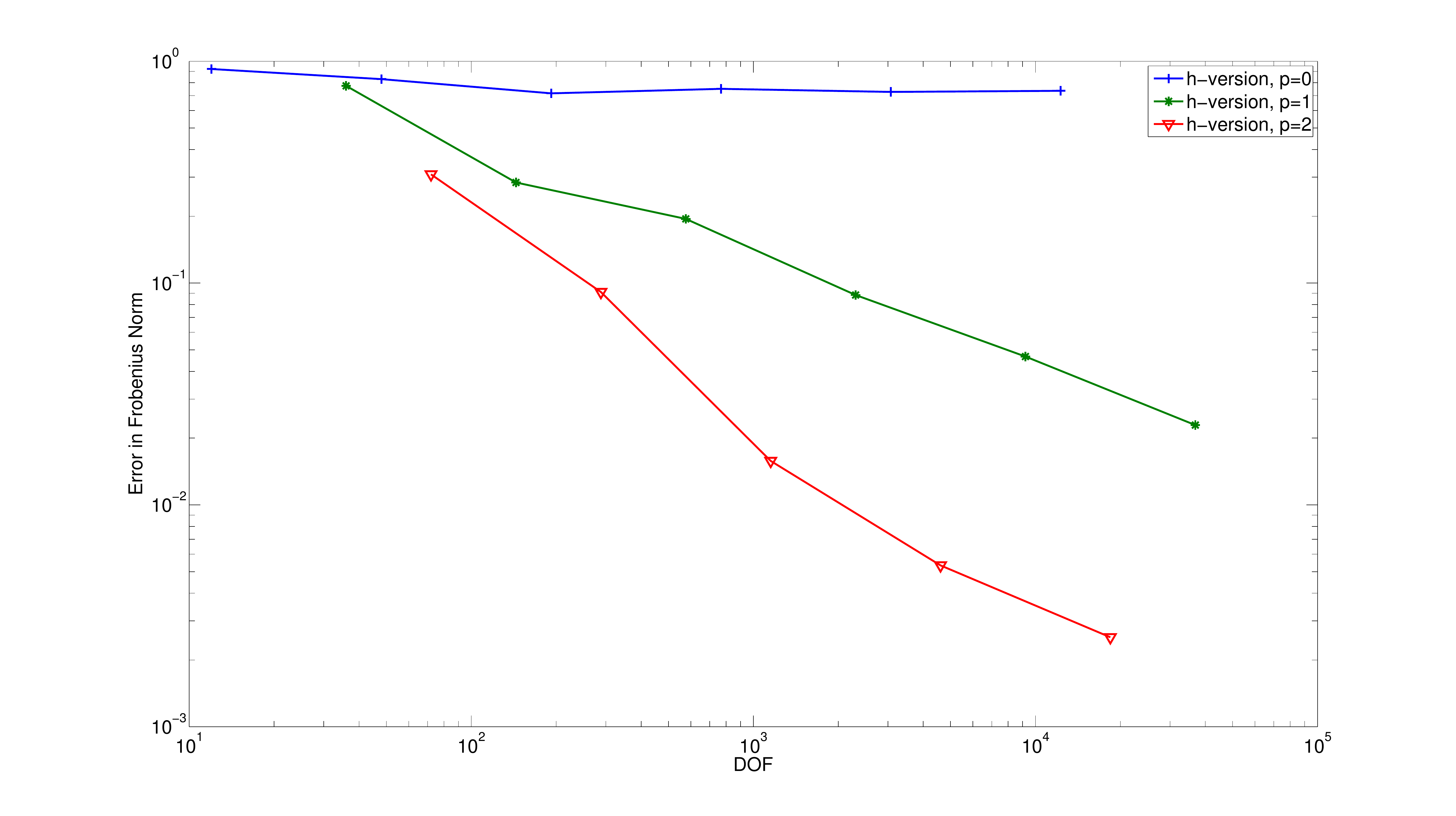}
	\caption{Error of the Hessian approximation in 3D for a point on the boundary surface}
	\label{fig:hessian3D}
\end{figure}

\begin{example}
Let $\Omega=[-1,1]^3$ be the domain and $g$ corresponding to the exact solution $u(x)=\frac{1}{\left\|x\right\|}$ with Hessian $H(x)= \frac{3}{\left\|x\right\|^5} \left(
\begin{array}[pos]{c c c}
	x_1^2 & x_1x_2 & x_1x_3 \\ x_1x_2 & x_2^2 & x_2x_3 \\ x_1x_3 & x_2x_3 & x_3^2
\end{array}\right)-\frac{1}{\left\|x\right\|^3} I$.
Figure~\ref{fig:hessian3D} shows the error of the Hessian in $x=(1,\frac{1}{3},\frac{1}{3})$. For $p \geq 2$ we observe good convergence of the Hessian approximation.
\end{example}

\section{Numerical Experiments}
\label{numex}
For the numerical experiments we set $\varphi: \mathbb S^2 \rightarrow \mathbb R^3$ to be $\varphi(x)=1.1 x$. This means that the sought surface is a sphere of radius $1.1$ with gravity potential $W_{meas}=\frac{1}{1.1}$ and gravity vector $G_{meas}=-\frac{1}{1.1^2} \frac{x}{|x|}$, both defined on $\mathbb S^2$.
The initial approximation $\varphi_0$ is the unit sphere $\mathbb S^2$. Therefore, $W_0=1$, $G_0=-\frac{x}{|x|}$ and $h_0=\frac{x}{2}$.

The sphere $\mathbb S^2$ is approximated by a regular, quasi-uniform mesh consisting of triangles such that the nodes of each triangle lie on $\mathbb S^2$.
More precisely, the mesh defines an icosahedron which is generated by \textit{maiprogs} \cite{01technicalmanual}. This mesh yields a domain approximation error and is kept fixed
for the entire Nash-H\"ormander algorithm. The main advantage is that only the coordinates of the nodes have to be updated and not the entire mesh itself.
This corresponds to a continuous, piecewise linear representation of $\varphi_m$, the new surface at the \textit{m}-th update of the algorithm.

The polynomial degree on each triangle is $p=2$, and $h_m$ in \eqref{eq:h_equation} is represented by a discontinuous piecewise constant function via interpolation in the midpoints of each triangle. Furthermore, $G_m$ is the linear interpolation in the nodes of the BEM approximation of $g|_{\varphi_m(\mathbb S^2)}$, obtained via $\bar{v}_m$ from equation \eqref{Dirichlet}.
Furthermore, since $h_m$ and the normal on each triangle $\mathcal T_n$ are piecewise constant, the jump contributions can be easily computed analytically. Furthermore, the operator $K'(h)$
in \eqref{definitionofS} can be computed semi-analytically by computing the action of the dual operator $K(h)$ on the test functions analytically \cite {MaischakTech} and performing an hp-composite Gaussian quadrature  \cite{Schwab} for the outer integration.

Since the boundary element space is the same for both the linear Molodensky problem \eqref{eq:linNLmol} and the auxiliary Dirichlet problem \eqref{Dirichlet}, the same single layer potential matrix is used in both \eqref{eq:discreteform1} and \eqref{eq:discreteform2}.
However, the computation of the right hand side $w_m$ for the Dirichlet problem is very CPU time consuming if a direct computation by means of \eqref{eq:SingleAnsatz}  is used. Since the ansatz and test functions live on varying surfaces, the computation
of one summand in \eqref{eq:SingleAnsatz} is as expensive as a semi-analytic computation of a single layer potential matrix. In particular, the computational time for the right hand side increases linearly with the number of iterations.

Since $\varphi_m$ is piecewise linear, the Gauss quadrature nodes $x$ for the outer integration are always mapped to exactly the same point on $\varphi_i(\mathbb S^2)$ under the mapping $\varphi_i \circ \varphi_m^{-1} (x)$
for each iteration step $m$. Therefore, if enough memory is available, $V_i\mu_i(\varphi_i \circ \varphi_m^{-1} (x))$ needs only be computed once and is stored for all the following iterations, keeping the
computational time for the right hand side  \eqref{eq:discreteform2} constant for all iterations $m$. This optimization together with the following  parallelization of the code  leads to a tremendous reduction of computing time.

With the solution of the Dirichlet problem \eqref{eq:discreteform2} at hand, the update of the surface in the nodes can be performed as defined in equation \eqref{corrphi} and with $g, \nabla g$ computed as in Section \ref{sec:3}.
The computation of one iteration is very CPU time consuming and therefore, parallelization of the code is crucial. Without parallelization and optimization of the code we need
$4+2m$ hours for the $m$-th iteration. However, with parallelization and optimization we need only $20$ minutes for each of the $m$ iterations for $N=2$~-icosahedron refinements corresponding to $320$ triangles whereas
we need $3$ hours for each of the $m$ iterations for $N=3$~-icosahedron refinements corresponding to $1280$ triangles. The numerical experiments were carried out on a cluster  with 5 nodes \`{a} 8 cores  with  2.93Ghz and 48GB memory,
where each core uses two Intel Nehalem X5570 processors.

In the following three different numerical experiments are presented. The first and the second experiment use  the classical Nash-H\"ormander algorithm with and without smoother as described in Section \ref{sec:2}. For the third experiment, the following restarted
algorithm with smoother is used.

\begin{alg} \label{alg:nashhormsmResart}
\vspace*{10mm} \hrule
\vspace*{0.5mm} (Nash-H\"ormander algorithm with restart) \vspace*{2mm}
\hrule
\begin{enumerate}
	\item For given measured data $W, G$ and $k\in \mathbb{R}$, choose $W_0, G_0,h_0,\varphi_0, \theta_0 \gg 1, \kappa \gg 1$
  \item Compute $W_k$, $G_k$, $\varphi_k$ in Algorithm~\ref{alg:nashhormsm}
	\item Stop if $\|| G_k-G\|| +\|| W_k-W\|| < tol$
	\item Else set	$W_0=W_k$, $G_0=G_k$, $h_0=h_k$ $\varphi_0=\varphi_k$, choose $\kappa, \theta_0$, and go to 2.
\end{enumerate}
\hrule
\end{alg}

\begin{figure}[tbp]
	\centering
	\includegraphics[keepaspectratio,width=120mm]{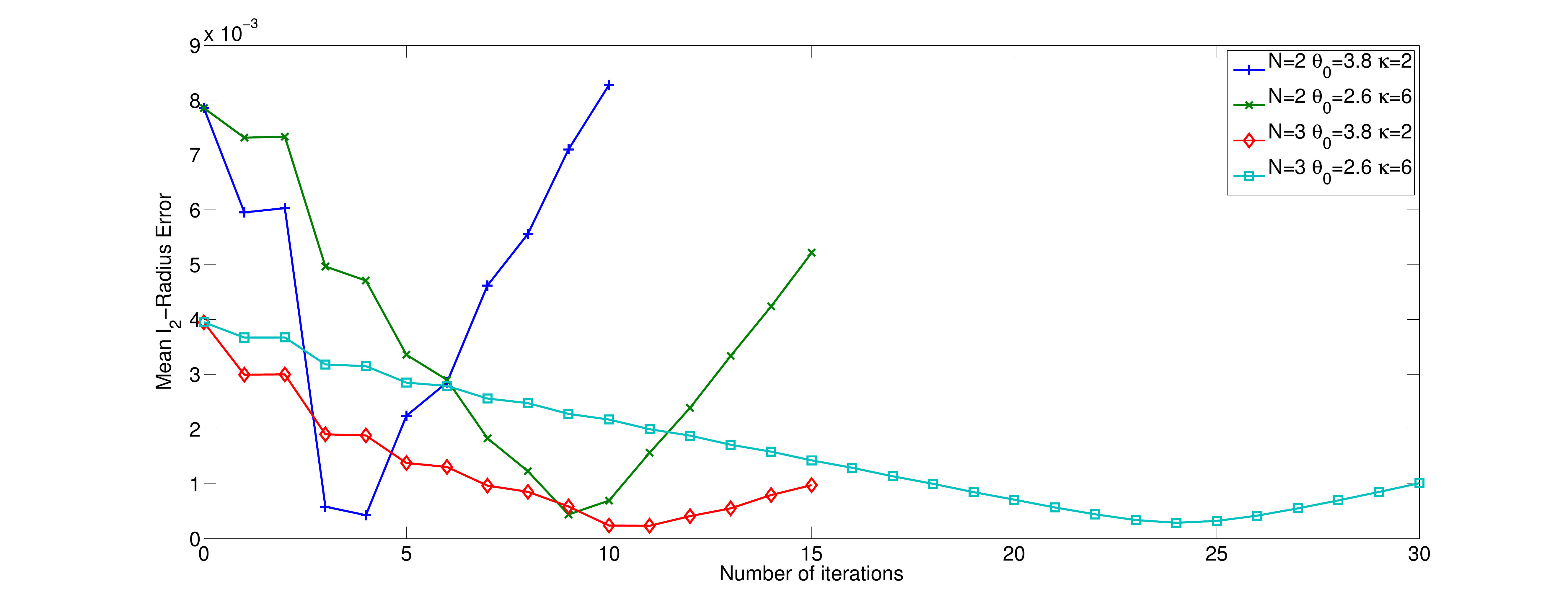}
	\caption{Mean Radius-Error in the $l_2$-norm with smoother}
	\label{fig:radiusl2smoother}
\end{figure}

Since the sought surface is also a sphere, we can expect that the sequences of computed surfaces are slightly perturbed spheres as well.
The perturbation should be a direct result of the domain approximation, different discretization errors and rounding errors. \\
Figure \ref{fig:radiusl2smoother} displays the mean $l_2$ error of the radius defined as

\begin{equation*}
\|e_r\|=\frac{1}{nr. nodes}\big[ \sum_{i=1}^{nr. nodes} (\|nodes(i)\|_2 - 1.1)^2 \big]^{1/2}
\end{equation*}
versus the number of iterations of Algorithm~\ref{alg:nashhormsm}. Whereas the algorithm itself takes care of the linearization error introduced by the linearization of the Molodensky problem, the algorithm does not tread the propagation of the discretization errors. Therefore, from a certain iteration step onwards the propagation of the spacial discretization error, for solving the linearized Molodensky problem, the auxiliary Dirichlet problem and computing the Hessian approximately, becomes dominating. Refining the mesh reduces the error and yields only a mild increase of the error for large iteration numbers. We have performed several numerical experiments with different parameters $\theta_0$ and $\kappa$. Firstly, if the amount of data smoothing is too small, the algorithm is unstable as expected. Secondly, if the amount of data smoothing is too large, then the essential information in the right hand side in the linearized Molodensky problem is lost in the first steps and in combination with the numerical errors convergence is lost. Also if the amount of smoothing does not decay sufficiently fast, the right hand side in the linearized Molodensky problem is close to machine precision leading to an ill-conditioned Hessian.

Here, we shortly comment how to perform smoothing with the heat kernel as it is used in the above example.
To smooth an arbitrary function  $F$, the heat equation with the Laplace-Beltrami operator is solved, where $F$ is the initial data.
 \begin{alignat*}{2}
\frac{\partial}{\partial t} u(x,t) - \varDelta u(x,t)&=0& &\,\,\,\mbox{in} \quad \varphi_m(\mathbb S^2) \times (0, \infty) ,\\
u(x,0) &=F(x)&   &\,\,\,\mbox{in} \quad  \varphi_m(\mathbb S^2).
\end{alignat*}
The unique solution of this problem is given by
\begin{equation}\label{heatsmoothu}
u(x,t)=\sum_{j=0}^{\infty} e^{-\lambda_j t} \langle F, \psi_j\rangle \psi_j(x).
\end{equation}
At $t=0$ we have
\begin{equation*}
u(x,0)=\sum_{j=0}^\infty \beta_j \psi_j(x)=F(x),
\end{equation*}
where $\beta_j$ are the Fourier coefficients $\langle F, \psi_j \rangle$. Here  $0=\lambda_0 < \lambda_1\leq \lambda_2 \leq \dots$ are the eigenvalues and $\psi_0, \psi_1, \psi_2, \dots $ the corresponding eigenfunctions
for the Laplace-Beltrami operator $\Delta$, i.e.~there holds
 \begin{equation} \label{laplacebeltramieigen}
  \Delta \psi_j =-\lambda_j \psi_j.
 \end{equation}
The eigenfunctions $\psi_j$ form an orthonormal basis in $L^2(\varphi(\mathbb S^2))$.

Having the discretized surface, \eqref{laplacebeltramieigen} can be solved approximately using the FEM method with continuous piecewise linear polynomials leading
to the generalized eigenvalue problem with the stiffness matrix $C$ and the mass matrix $A$ of the Laplace-Beltrami operator $\Delta$
\begin{equation}\label{eigeneq}
C \psi_h = \lambda_h A \psi_h
\end{equation}
where $\psi_h$ denotes the unknown $L^2$-orthonormal eigenfunction, evaluated at the mesh vertices. With $\psi_h$ solving \eqref{eigeneq} the heat kernel can be approximated by
 \begin{equation*}
e^{-t \Delta}(x,y)=\sum_{j=0}^M e^{-t \lambda_{j,h}} \psi_{j,h}(x) \psi_{j,h}(y)
\end{equation*}
where $M$ must be sufficiently large. Once we obtained the components  $\psi_{j,h}$ of the eigenfunctions $\psi_h$, we compute the Fourier coefficients $\beta_{j,h}$ as presented in \cite[Eqn. (10)]{Seo:2010:HKS:1926877.1926943}. Therewith,
\begin{equation} \label{heatuhxt}
u_h(x,t)= \sum_{j=0}^M e^{-\lambda_{j,h}t} \beta_{j,h} \psi_{j,h}(x).
\end{equation}
For our numerical experiments $F$ is always of the structure $G_{meas}\!-G_m\!+\!\! \sum_{j=0}^{m-1}\!\!\triangle_j \dot {\widetilde G}_j$ (see \eqref{greccurencetosmooth}).  We use $u_h(x, \frac{1}{{}\,\,\theta_m})$  where $t=\frac{1}{{}\,\,\theta_m}$ in \eqref{heatuhxt} as the smoothed $F$,
where $\theta_m$ is computed by \eqref{thetatocompute}.

\begin{figure}[tbp]
	\centering
	\includegraphics[keepaspectratio,width=120mm]{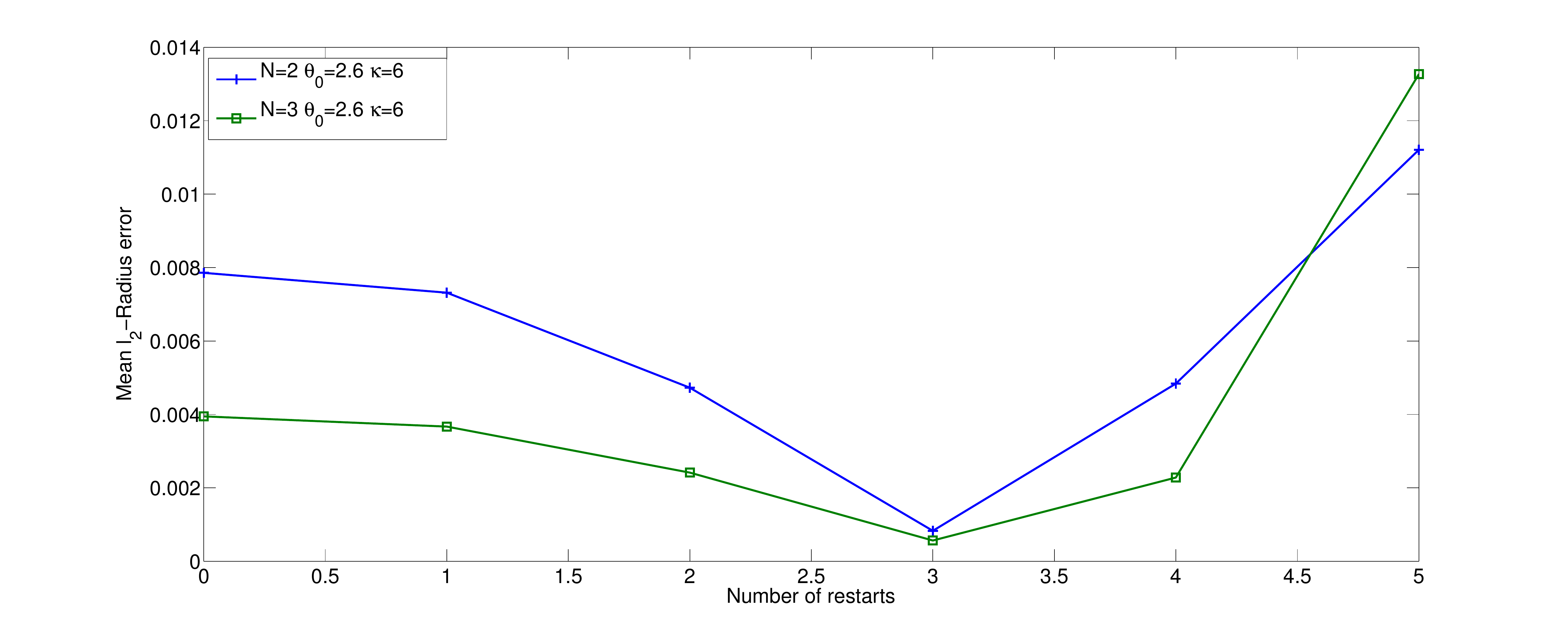}
	\caption{Mean Radius-Error  in the $l_2$-norm with smoother and restart}
	\label{fig:radiusrestart}
\end{figure}

Figure \ref{fig:radiusrestart} displays the mean $l_2$ error of the radius versus the number of restarts for the restarted algorithm presented in this section with smoother. In the extreme case, in which the algorithm is restarted after each iteration, we still observe the same structural behavior as for the other experiments. In particular, from the third restart onwards the discretization error
propagation becomes dominating again. Again the error can be reduced by refining the mesh.

\begin{table}[h!]
\begin{center}
\begin{tabular}{|r|r|l|l|l|}
\hline
Iter  & DOF & $|(u_N(\boldsymbol q)-u(\boldsymbol q))|$ & EOC \\
\hline
0&  120    &\hspace{0.4cm}1.04164e+03&  \\
  &  480    &\hspace{0.4cm}3.94212e+02& 0.70  \\
  &1920    & \hspace{0.4cm}1.04676e+02 & 0.96 \\
  &7680    & \hspace{0.4cm}27.79515& 0.96 \\
 \hline
1&  120    &\hspace{0.4cm}5.82518e+03&  \\
   &  480    &\hspace{0.4cm}3.66919e+02& 1.99  \\
   &1920    & \hspace{0.4cm}1.05925e+02&  0.90\\
   &7680    & \hspace{0.4cm}30.57947& 0.90\\
\hline
2&  120    &\hspace{0.4cm}2.96617e+03&  \\
  &  480    &\hspace{0.4cm}1.01239e+03& 0.77 \\
  &1920    & \hspace{0.4cm}2.72407e+02&  0.95\\
  &7680    & \hspace{0.4cm}73.29735& 0.95\\
\hline
 \end{tabular}
\caption{Pointwise Errors for the linearized Molodensky problem with smoother}
\label{tab:errorspointwise}
\end{center}
\end{table}



\begin{figure}[h!]
	\centering
	 \subfigure[N=2, 320 triangles, $r=1.007$]{
		\includegraphics[width=52.5mm, keepaspectratio]{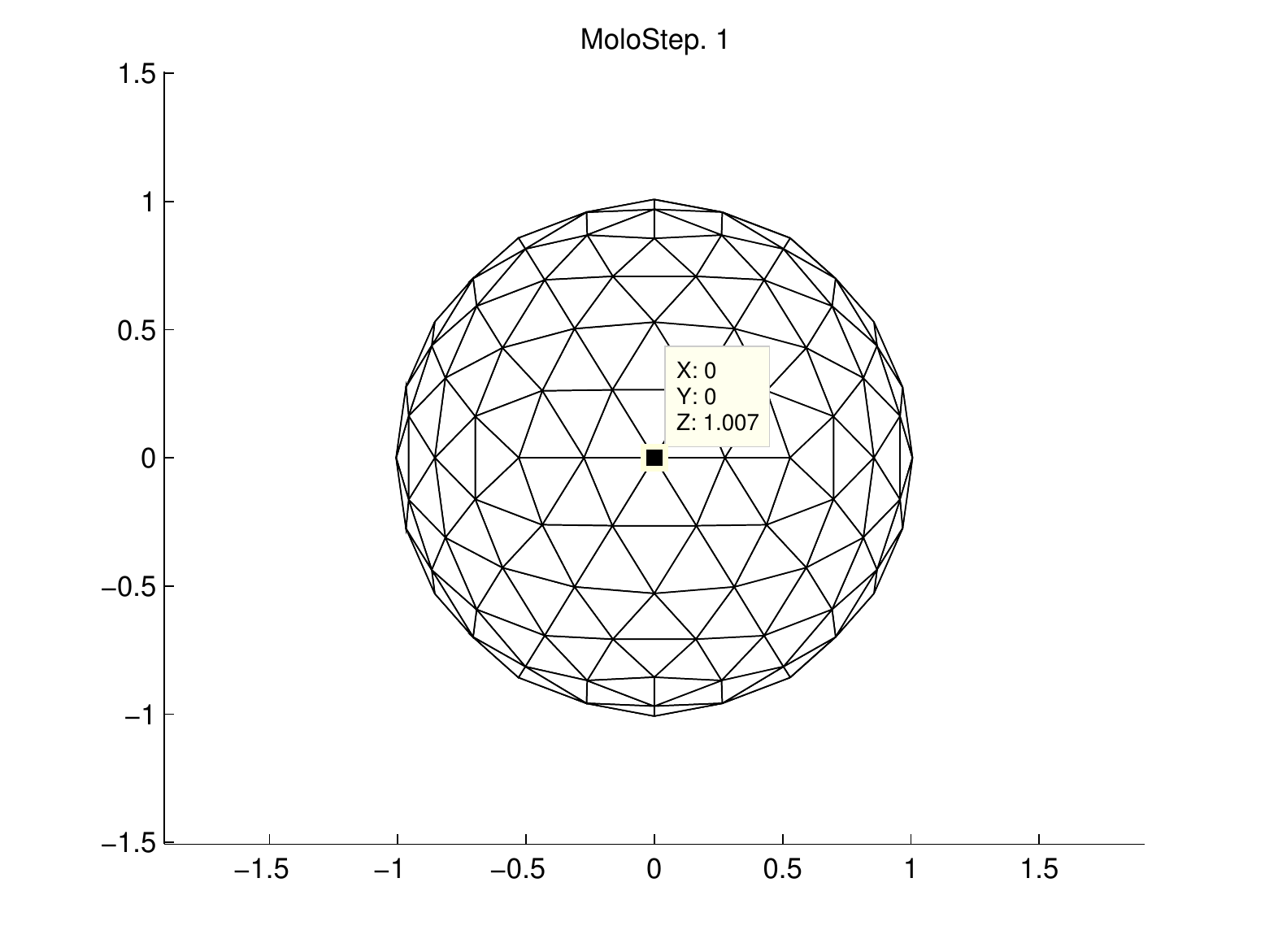}} \
	 \subfigure[N=3, 1280 triangles, $r=1.007$]{
		\includegraphics[width=52.5mm, keepaspectratio]{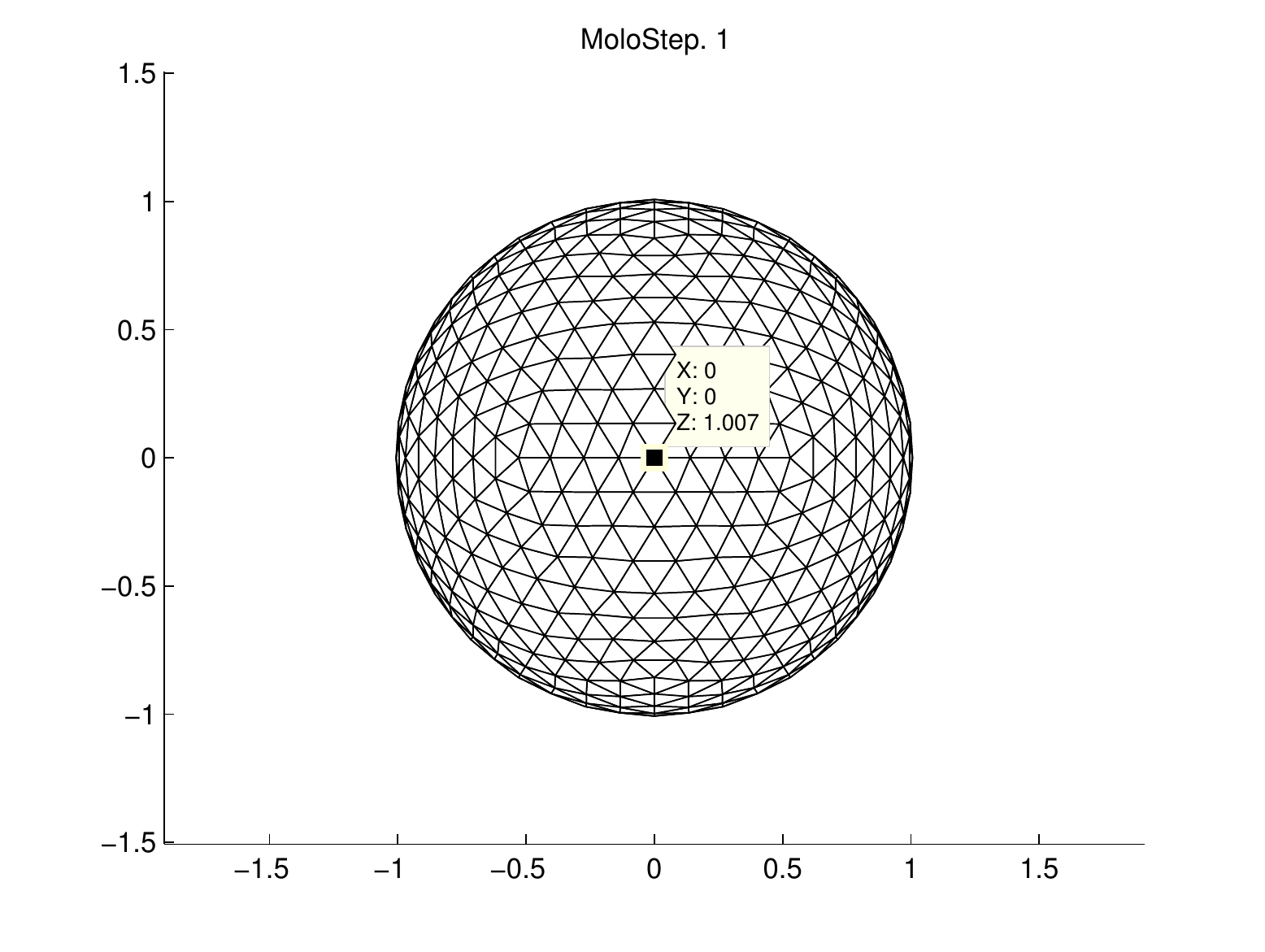}}\\
	 \subfigure[N=2, 320 triangles, $r=1.057$]{
		\includegraphics[width=52.5mm, keepaspectratio]{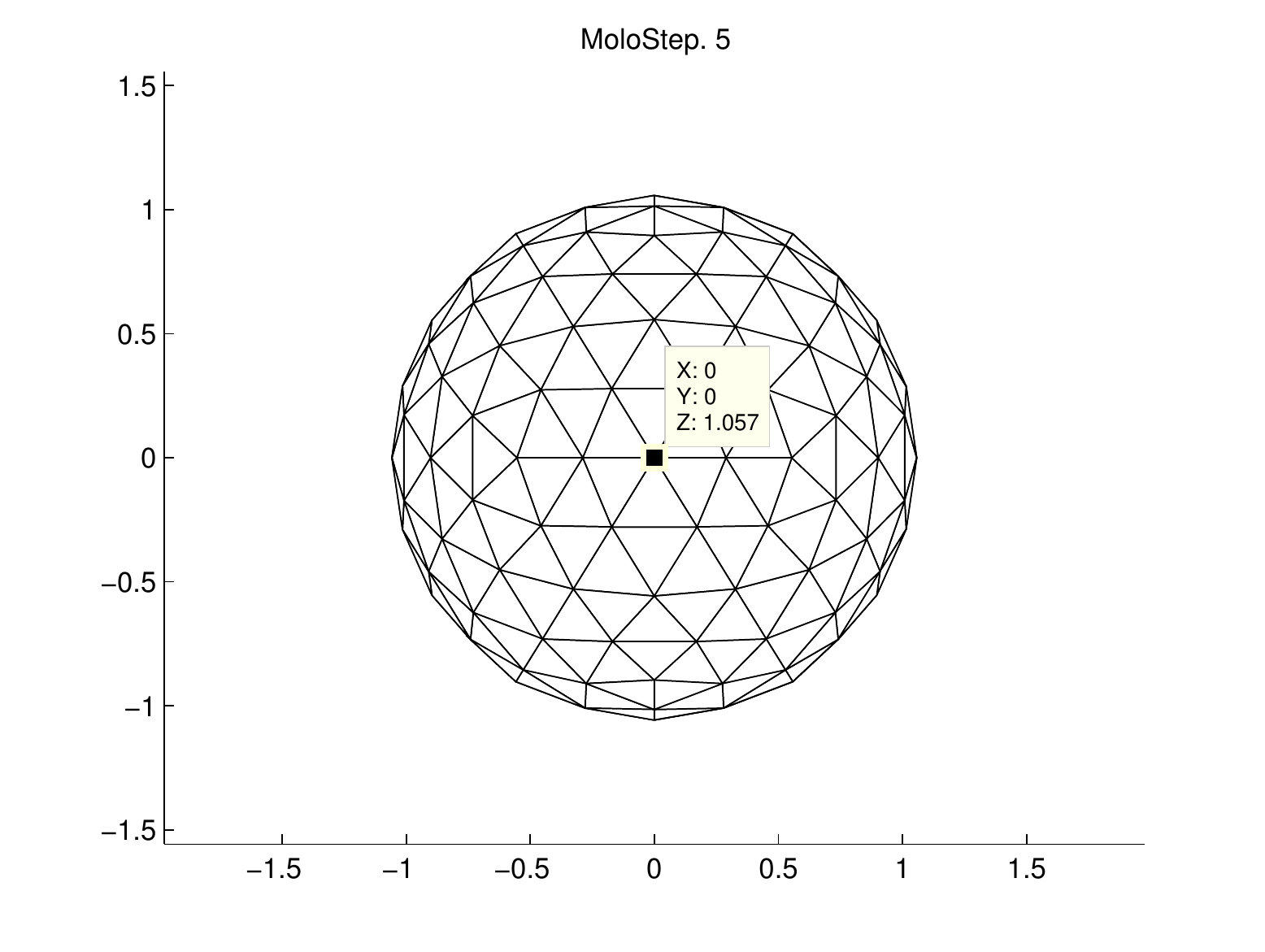}} \
	 \subfigure[N=3, 1280 triangles, $r=1.045$]{
		\includegraphics[width=52.5mm, keepaspectratio]{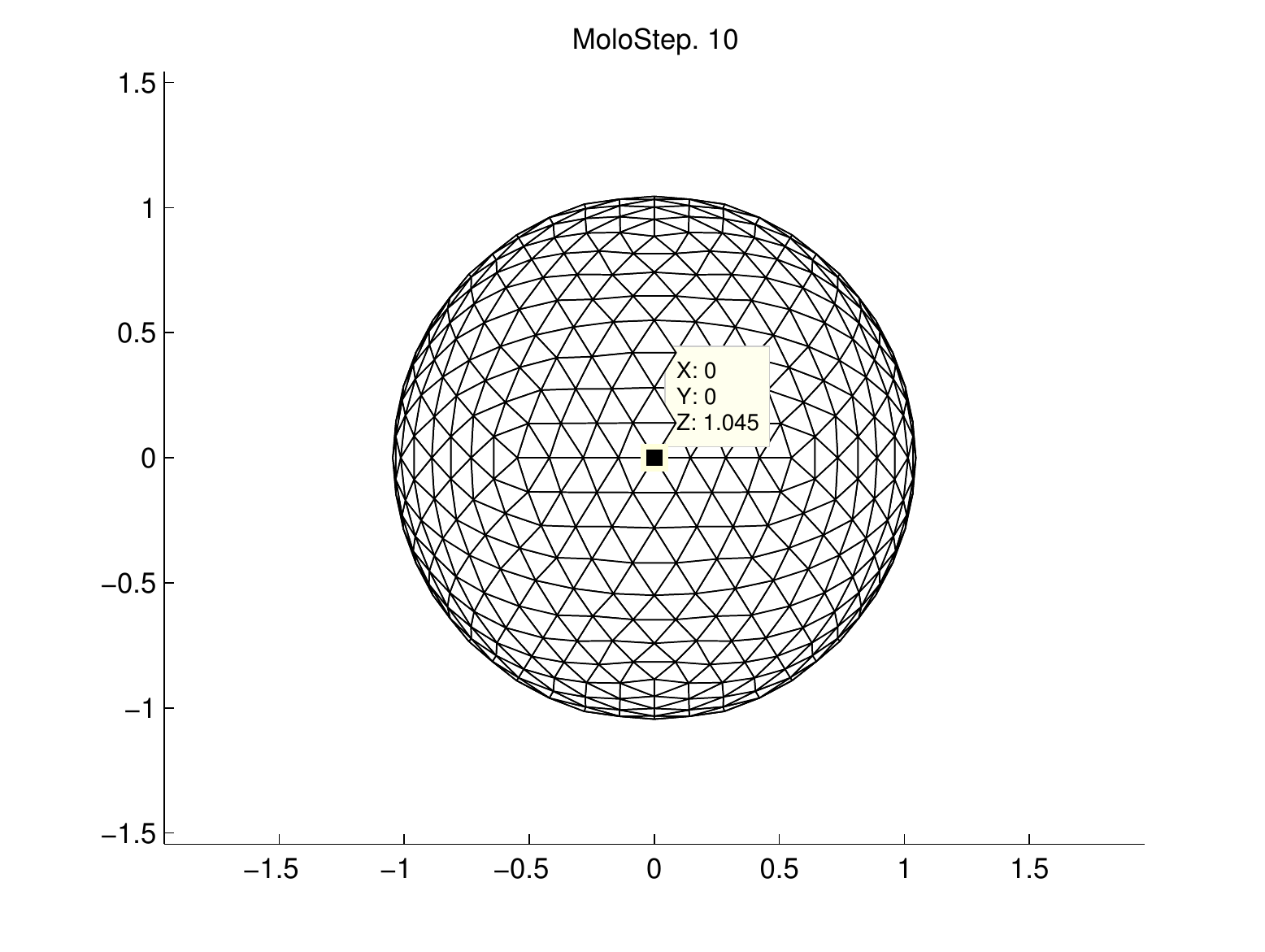}}\\
 \subfigure[N=2, 320 triangles, $r=1.107$]{
		\includegraphics[width=52.5mm, keepaspectratio]{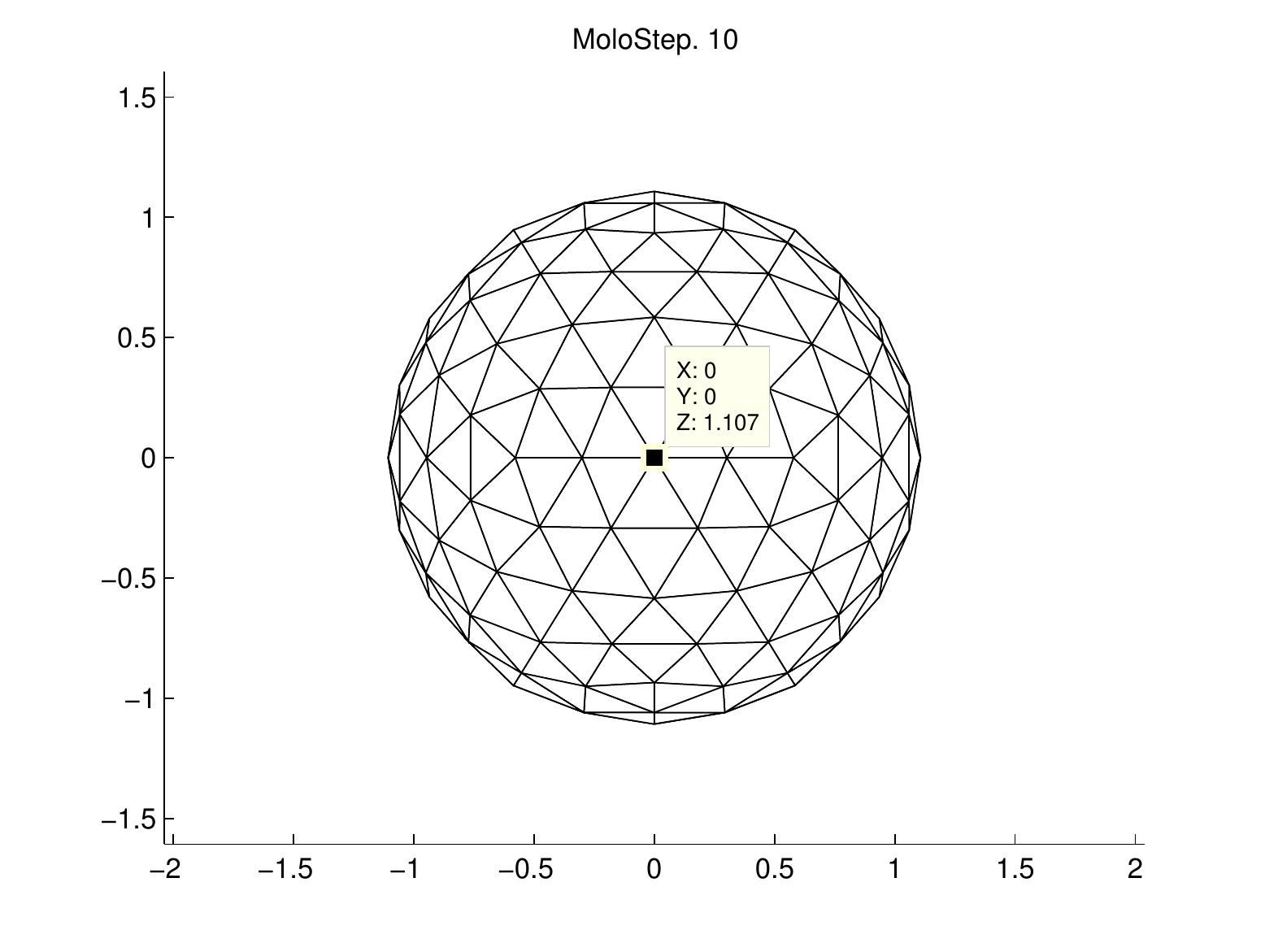}} \
	 \subfigure[N=3, 1280 triangles, $r=1.104$]{
		\includegraphics[width=52.5mm, keepaspectratio]{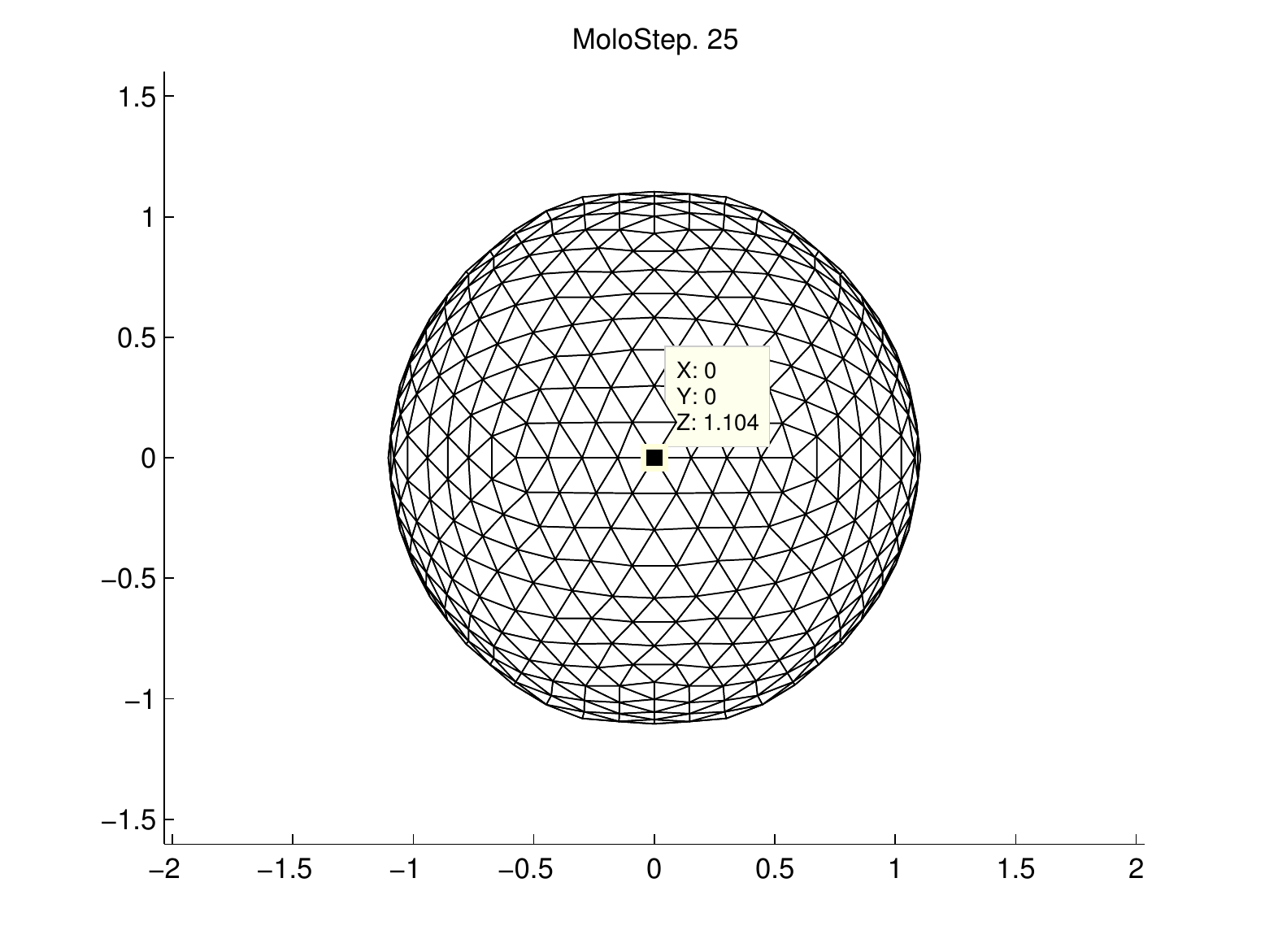}}		
		\caption{$N=2$ and $N=3$ icosahedron refinements with smoother, $\theta_0=2.6, \kappa=6$}
	\label{fig:figuressm2}
\end{figure}

Figure \ref{fig:figuressm2} displays the sequence of obtained spheres.


For large scale applications one has to use standard reduction methods for BEM to reduce the computational complexity, e.g.~$H$-matrices and/or multipole expansion or wavelet compression techniques. All these techniques can be applied to speed up our algorithm and allow computations on finer meshes.

\section*{Acknowledgments}
This work was supported by the cluster of excellence QUEST, the Danish National Research Foundation (DNRF) through the Centre for Symmetry and Deformation and the Danish Science Foundation (FNU) through research grant 10-082866. 

\bibliographystyle{abbrv}
\bibliography{molodensky_numeric_paper}

\end{document}